# ON NONPARAMETRIC AND SEMIPARAMETRIC TESTING FOR MULTIVARIATE LINEAR TIME SERIES

BY YOSHIHIRO YAJIMA[1] AND YASUMASA MATSUDA[2]

*University of Tokyo and Tohoku University*

We formulate nonparametric and semiparametric hypothesis testing of multivariate stationary linear time series in a unified fashion and propose new test statistics based on estimators of the spectral density matrix. The limiting distributions of these test statistics under null hypotheses are always normal distributions, and they can be implemented easily for practical use. If null hypotheses are false, as the sample size goes to infinity, they diverge to infinity and consequently are consistent tests for any alternative. The approach can be applied to various null hypotheses such as the independence between the component series, the equality of the autocovariance functions or the autocorrelation functions of the component series, the separability of the covariance matrix function and the time reversibility. Furthermore, a null hypothesis with a nonlinear constraint like the conditional independence between the two series can be tested in the same way.

**1. Introduction.** One of the main purposes of multivariate stationary time series analysis is to elucidate intrinsic relationships between different component series. Frequently, these relationships can be expressed in terms of specific constraints imposed on the spectral density matrix. For instance, the spectral density matrix of a separable time series is a product of the contemporaneous covariance matrix of the component series and the scalar spectral density function, which is common to them (see, e.g., Haslett and Raftery [14], Martin [23], Cressie [6], Guyon [13], Shitan and Brockwell [36] and Matsuda and Yajima [25]). If the underlying time series is Gaussian, the independence between the component series is equivalent, so that the

Received October 2006; revised February 2008.

[1]Supported in part by JSPS Grant A(1) 15200021.

[2]Supported in part by JSPS Grant C(2)14580345 A(1)15200021.

spectral density matrix is diagonal for all frequencies (see, e.g., Wahba [40]). The conditional independence of two components, given the others, is equivalent, so that the corresponding partial spectral coherence is identical to zero (see, e.g., Dahlhaus [7]). The time reversibility is characterized by the spectral density matrix being real-valued (see, e.g., Chan, Ho and Tong [3]). Furthermore, the equality of the autocovariance functions or the autocorrelation functions of the component series implies that their spectral density functions are equal to each other with an appropriate scale-shift for the latter case (see, e.g., Coates and Diggle [5] and Diggle and Fisher [10]).

In this paper, we formulate nonparametric and semiparametric hypothesis testing on multivariate linear time series in a unified fashion and propose new test statistics based on estimators of the spectral density matrix. First, we construct two estimators of the spectral density matrix, of which the first is always consistent and the latter is consistent only when the null hypothesis is true, and, next, we introduce a function to measure the discrepancy between these estimators. Then, the discrepancy is asymptotically standard normally distributed under the null hypothesis as the sample size goes to infinity, whereas, if the null hypothesis is false, it diverges to infinity. Related ideas are applied for discriminant analysis of time series by Kakizawa, Shumway and Taniguchi [18] and for nonparametric testing on univariate time series by Taniguchi and Kondo [38] and multivariate time series by Taniguchi, Puri and Kondo [39], respectively (see Taniguchi and Kakizawa [37] for a comprehensive exposition). Hong and White [17] considered specification testing on nonparametric and semiparametric regression models in a similar way.

On the other hand, for parametric models, Paparoditis [29] and [30] and Delgado, Hidalgo and Velasco [8] considered test statistics based on another discrepancy measure between the hypothesized spectral density function (or matrix) and their estimators.

The advantages of our test statistics are as follows. First, our test statisctics can be applied to both nonparametric and semiparametric hypotheses on multivariate linear time series in a unified way. They are robust against wrong decisions caused by misspecification, which parametric approaches often suffer.

Second, under the null hypothesis, the limiting distributions of the test statistics are always normal distributions. Moreover, since they are scale invariant, that is, independent of the unit of measurement for the observations, the expectations and variances of the limiting distributions can be expressed in a simple form, interestingly, some of them are known constants in the preceding examples. Hence, they can be implemented easily for practical use. In contrast, the limiting distributions of some test statistics, which have already been proposed for specific hypotheses mentioned above, are unknown or have rather complicated forms, particularly when a nonlinear constraint is imposed on the null hypothesis, which makes it difficult to show the appropriate validity of their significance tests (see, e.g., Dahlhaus [7]).

Finally, if the statistics are normalized to be standard normal asymptotically under the null hypothesis, then, under any alternative, they diverge to infinity and consequently are consistent tests. However, a trade-off for this advantage is that our test can only detect local alternatives of $O(n^{-(1+\beta)/4})$ with $1/2 < \beta < 3/4$ being slightly slower than $n^{-1/2}$, which is shared with the test for regression models proposed by Hong and White [17].

The paper is organized as follows. We introduce the mathematical formulation of the null hypothesis and the test statistics in Section 2. In Section 3, we derive the theorems on limiting properties of the test statistics. In Section 4, we apply the theorems to the examples mentioned above. Some simulation results are shown in Section 5. Section 6 is devoted to prove the theorems. The technical lemmas are listed in Section 7.

The mathematical details on the lemmas, the theorems, the examples and the computational simulations are available from the authors upon request.

## 2. Null hypotheses and test statistics.

Let $\mathbf{Z}_t = (Z_{1t}, \dots, Z_{rt})'$ be an $r$-dimensional zero-mean stationary Gaussian time series possessing a spectral density matrix $f(\lambda)$. $f(\lambda)$ is defined as periodic with the period $2\pi$ for $\lambda \notin (-\pi, \pi]$.

Throughout this paper, $A_{ab}$ and $A^{ab}$ are generic symbols for the $(a, b)$th element of matrices $A$ and $A^{-1}$, respectively. $A'$ is the transposed matrix of $A$.

Let $g(\theta, y) = (g_{ab}(\theta, y))$ be an $r \times r$ matrix-valued function, where $y = (y_{\alpha\beta})$ is an $r \times r$ matrix and $\theta = (\theta_1, \dots, \theta_v)'$ is a $v$-dimensional vector of parameters. We assume that, under the null hypothesis, $f(\lambda)$ satisfies the equation

$$(1) \qquad f(\lambda) = g(\theta, f(\lambda))$$

for all $\lambda \in (-\pi, \pi]$.

If $g_{ab}(\theta, y) = y_{ab}$, for all $a$ and $b$, it implies that there is no constraint imposed on $f(\lambda)$. Otherwise, (1) introduces some relationship between the components of $f(\lambda)$. The hypothesis is semiparametric for $v > 0$ and is nonparametric for $v = 0$, respectively.

Now, we define the test statistics. Given $n$ observations $\mathbf{Z}_1, \dots, \mathbf{Z}_n$, introduce the discrete Fourier transform $W_a(\lambda) = \frac{1}{\sqrt{2\pi n}} \sum_{t=1}^n Z_{at} \exp(it\lambda)$ for $a = 1, \dots, r$. The cross periodogram of $Z_{at}$ and $Z_{bt}$ is

$$I_{Z,ab}(\lambda) = W_a(\lambda)\overline{W_b(\lambda)}, \qquad a, b = 1, \dots, r,$$

and the periodogram matrix is $I_Z(\lambda) = (I_{Z,ab}(\lambda))$. $I_Z(\lambda)$ is defined as periodic with the period $2\pi$ for $\lambda \notin (-\pi, \pi]$. Denote the Fourier frequency $\frac{2\pi j}{n}$ as $\lambda_j$ for $j = 0, \pm, 1, \ldots$. For notational simplicity, $W_a(\lambda_j)$, $I_{Z,ab}(\lambda_j)$ and $I_Z(\lambda_j)$ are denoted by $W_{a,j}$, $I_{Z,ab,j}$ and $I_{Z,j}$, respectively. Actually, they depend on the sample size $n$ but we suppress it.

Put $f_t = f(\lambda_t)$, and define the unrestricted estimator of $f_t$ by the smoothed periodogram matrix $\hat{f}_{U,t} = (\hat{f}_{U,ab,t})$, where

$$(2) \qquad \hat{f}_{U,t} = \frac{1}{w^*} \sum_{j=-m/2}^{m/2} w_j I_{Z,t+j}, \qquad t = 1, \ldots, [n/2],$$

and $w_j$ $(j = -m/2, \ldots, m/2)$ is a weight sequence with $w^* = \sum_{j=-m/2}^{m/2} w_j$ and $[x]$ is the integer part of $x$.

Under the null hypothesis, we define the estimator $\hat{f}_{R,t} = (\hat{f}_{R,ab,t})$ by

$$\hat{f}_{R,t} = g(\hat{\theta}_n, \hat{f}_{U,t}), \qquad t = 1, \ldots, [n/2],$$

where $\hat{\theta}_n = (\hat{\theta}_{1n}, \ldots, \hat{\theta}_{vn})'$ is an estimator of $\theta$.

We can expect that, under the null hypothesis (1), $\hat{f}_{U,t}\hat{f}_{R,t}^{-1}$ is close to $I_r$, the $r \times r$ identity matrix, whereas they are far from each other if (1) is false. Hence, we introduce a function that measures the discrepancy between $\hat{f}_{U,t}\hat{f}_{R,t}^{-1}$ and $I_r$. Let $K(A)$ be a nonnegative function defined on all $r \times r$ complex matrices that are similar to a positive definite matrix, and let it be zero if and only if $A = I_r$. Note that $\hat{f}_{U,t}\hat{f}_{R,t}^{-1}$ is similar to $\hat{f}_{R,t}^{-1/2}\hat{f}_{U,t}\hat{f}_{R,t}^{-1/2}$.

Then, we introduce the statistic

$$T_n = \sum_{t=1}^{[n/2]} K(M_t),$$

where $M_t = (m_{ab,t}) = \hat{f}_{U,t}\hat{f}_{R,t}^{-1}$.

Now, we give a few candidates for $K(A)$ (see [18]). The first one is Kullback–Leibler (KL) discrimination information defined by

$$K_I(A) = \operatorname{tr}(A) - \log \det(A) - r$$

(see Kullback and Leibler [20] and Kullback [19]). Next, the *J divergence*, which is a symmetric version of the KL discrimination information, is

$$K_J(A) = K_I(A) + K_I(A^{-1})$$

(see Kullback and Leibler [20] and Kullback [19]). Finally, using the Chernoff information (see Chernoff [4] and Renyi [33]), Parzen [31] proposed

$$K_\alpha(A) = \log \det(\alpha A + (1-\alpha)I_r) - \alpha \log \det(A), \qquad 0 < \alpha < 1.$$

We shall show in the subsequent section that, under the null hypothesis,

$$\sqrt{\frac{m}{n}}\left(T_n - \frac{n}{m}\eta\right)/\sigma$$

is asymptotically standard normally distributed with some constants $\eta$ and $\sigma$. Consequently if $\hat{\eta}_n$ and $\hat{\sigma}_n$ are consistent estimators of $\eta$ and $\sigma$, respectively, and $\hat{\eta}_n - \eta = o_p(m^{1/2}/n^{1/2})$, then the test statistic

$$\hat{T}_n = \sqrt{\frac{m}{n}}\left(T_n - \frac{n}{m}\hat{\eta}_n\right)/\hat{\sigma}_n$$

is also asymptotically standard normally distributed. Actually, $\eta$ and $\sigma^2$ of some nonparametric examples mentioned in Section 1 are known constants, and they require no estimation procedure. On the other hand, $\hat{T}_n$ diverges to infinity if (1) is false and, hence, is a consistent test for any alternative.

**3. Theorems.** Before we proceed to show the main theorems, we introduce some assumptions and notation. First, we introduce the following assumptions:

(A1) $\mathbf{Z}_t$ is an $r$-dimensional zero-mean stationary Gaussian process;

(A2) $f(\lambda)$ is a positive definite matrix for all $\lambda \in (-\pi, \pi]$;

(A3) $f(\lambda)$ is twice continuously differentiable for all $\lambda \in (-\pi, \pi]$;

(A4) $m = O(n^\beta), \frac{1}{2} < \beta < \frac{3}{4}$, and the weight sequence $w_j, j = -m/2, \ldots, m/2$ is

$$w_j = u(j/m), \qquad j = -m/2, \ldots, m/2,$$

where $u(x)$ is a positive continuously differentiable even function on $[-1/2, 1/2]$;

(A5) $g(\theta, y)$ is three times continuously differentiable for $\theta$ and $y$;

(A6) $K(A)$ is four times continuously differentiable for $A$.

Next, put $\theta_0 = (\theta_{10}, \ldots, \theta_{v0})'$, $\breve{g}(\lambda) = (g(\theta_0, f(\lambda))$ and $\breve{g}_t = \breve{g}(\lambda_t)$ if the null hypothesis (1) is true. Finally, let $k_{ab,cd}^{(2)}$ be the second partial derivative $\partial^2 K(A)/\partial A_{ab}\,\partial A_{cd}$ evaluated at $A = I_r$. The third partial derivative $k_{ab,cd,ef}^{(3)}$ is defined similarly.

By the results in Appendix A.13 of Lütkepohl [22], for $K_I(A)$,

$$\kappa_{ab,cd}^{(2)} = \begin{cases} 1, & \text{if } a = d, b = c, \\ 0, & \text{otherwise.} \end{cases} \tag{3}$$

For $K_J(A)$ and $K_\alpha(A)$, their second partial derivatives are given by $2\kappa_{ab,cd}^{(2)}$ and $(\alpha - \alpha^2)\kappa_{ab,cd}^{(2)}$, respectively.

From now on, to avoid abuse of notation and give the theoretical results in a concise form, we consider $K_I(A)$ and assume (3), because, for other $K(A)$'s, the corresponding results are given in the same way without any modification.

Then, our main results are the following.

THEOREM 1. *Under assumptions* (A1)–(A6) *and ([3]), if the null hypothesis ([1]) is true and $\hat{\theta}_n - \theta_0 = O_p(n^{-1/2})$, then the limiting distribution of $(m/n)^{1/2}(T_n - (n/m)\eta)$ is $N(0, \sigma^2)$, where*

$$\eta = \sum_{\alpha,\beta,\gamma,\nu=1}^{r} \frac{C_u}{\pi} \int_0^\pi \mu_{\alpha\beta\gamma\nu}(\lambda) \breve{g}_{\alpha\nu}(\lambda) \breve{g}_{\gamma\beta}(\lambda)\, d\lambda,$$

$$\sigma^2 = \sum_{\alpha,\beta,\gamma,\nu=1}^{r} \sum_{\alpha',\beta',\gamma',\nu'=1}^{r} \frac{D_u}{\pi} \int_0^\pi \mu_{\alpha\beta\gamma\nu}(\lambda) \overline{\mu_{\alpha'\beta'\gamma'\nu'}(\lambda)}$$
$$\times (\breve{g}_{\alpha\alpha'}(\lambda) \breve{g}_{\beta'\beta}(\lambda) \breve{g}_{\gamma\gamma'}(\lambda) \breve{g}_{\nu'\nu}(\lambda)$$
$$+ \breve{g}_{\alpha\gamma'}(\lambda) \breve{g}_{\nu'\beta}(\lambda) \breve{g}_{\gamma\alpha'}(\lambda) \breve{g}_{\beta'\nu}(\lambda))\, d\lambda$$

*and*

$$\mu_{\alpha\beta\gamma\nu}(\lambda) = \frac{1}{2} \operatorname{tr}\left[ \breve{g}(\lambda) \frac{\partial \breve{g}^{-1}}{\partial y_{\alpha\beta}}(\lambda) \breve{g}(\lambda) \frac{\partial \breve{g}^{-1}}{\partial y_{\gamma\nu}}(\lambda) \right] + \frac{1}{2}\left[ \frac{\partial \breve{g}^{-1}}{\partial y_{\alpha\beta}}(\lambda) \right]_{\nu\gamma}$$
$$+ \frac{1}{2}\left[ \frac{\partial \breve{g}^{-1}}{\partial y_{\gamma\nu}}(\lambda) \right]_{\beta\alpha} + \frac{1}{2} \breve{g}^{\beta\gamma}(\lambda) \breve{g}^{\nu\alpha}(\lambda),$$

*and $C_u$ and $D_u$ are defined in Lemmas [9] and [10], respectively.*

Hence, if the null hypothesis ([1]) is true and $\hat{\eta}_n - \eta = o_p(m^{1/2}/n^{1/2})$ and $\hat{\sigma}_n^2 - \sigma^2 = o_p(1)$, then the limiting distribution of $\hat{T}_n$ is the standard normal distribution.

THEOREM 2. *Under assumptions* (A1)–(A6) *and ([3]), if the null hypothesis ([1]) is false, and $\hat{\theta}_n$, $\hat{\eta}_n$ and $\hat{\sigma}_n$ converge in probability to some constants $\theta^*$, $\eta^*$ and $\sigma^*$, respectively, so that $g(\theta^*, f(\lambda))$ is a positive definite matrix for all $\lambda$ in $(-\pi, \pi]$ and $\sigma^*$ is positive, then, for any sequence $\{C_n\}, C_n = o((nm)^{1/2})$,*

$$\lim_{n\to\infty} Pr[\hat{T}_n > C_n] = 1.$$

Theorems [1] and [2] assure that our test statistics are consistent for any alternative. Furthermore, if we put

$$C^* = \frac{1}{2\pi\sigma^*} \int_0^\pi K(f(\lambda)g(\theta^*, f(\lambda))^{-1})\, d\lambda,$$

the proof of Theorem [2] shows that, if the null hypothesis is false, $\hat{T}_n/\sqrt{nm}$ converges in probability to $C^*$ as $n \to \infty$.

Next consider asymptotic behavior of the test statistics under local alternatives. To state our theorem, we define a class of local alternatives

$$H_{an}: f_n(\lambda) = f(\lambda) + \frac{1}{(mn)^{1/4}} f^*(\lambda),$$

where $f(\lambda)$ satisfies (1) with $\theta = \theta_0$ and $f^*(\lambda)$ is a spectral density matrix. Then, we introduce the following assumptions, which specify the behaviors of $f^*(\lambda)$ and the parametric estimator $\hat{\theta}_n$:

(A7) $f^*(\lambda)$ is a positive definite matrix and twice continuously differentiable for $\lambda$ on $(-\pi, \pi]$;

(A8) Under the local alternative $H_{an}$, there exists a nonstochastic sequence $\{\theta_n^*\}$ such that $\hat{\theta}_n - \theta_n^* = O_p(n^{-1/2})$ and $\lim_{n \to \infty}(mn)^{1/4}(\theta_n^* - \theta_0) = \rho_0 = (\rho_{10}, \ldots, \rho_{v0})'$.

Before we proceed to the theorem, we introduce some notation. Set $\check{g}_n(\lambda) = g(\theta_n^*, f_n(\lambda))$ and $\check{g}_{t,n} = \check{g}_n(\lambda_t)$. $\eta_n$ and $\sigma_n^2$ are defined by substituting $\check{g}_n(\lambda)$ for $\check{g}(\lambda)$ in $\eta$ and $\sigma^2$ of Theorem 1, respectively. Then, we have the following result.

THEOREM 3. Under (A1)–(A8) and (3), if $H_{an}$ is true, $\hat{\eta}_n - \eta_n = o_p(m^{1/2}/n^{1/2})$ and $\hat{\sigma}_n^2 - \sigma_n^2 = o_p(1)$, then the limiting distribution of $\hat{T}_n$ is $N(\xi/\sigma, 1)$, where

$$\xi = \frac{1}{4\pi}\int_0^\pi \text{tr}\Bigg[\Bigg(\sum_{i=1}^v \rho_{i0}\bigg[\frac{\partial \check{g}(\lambda)}{\partial \theta_i}\bigg]\check{g}^{-1}(\lambda)$$
$$+ \sum_{\alpha,\beta=1}^r f_{\alpha\beta}^*(\lambda)\bigg[\frac{\partial \check{g}(\lambda)}{\partial y_{\alpha\beta}}\bigg]\check{g}^{-1}(\lambda) - f^*(\lambda)\check{g}^{-1}(\lambda)\Bigg)^2\Bigg]\,d\lambda.$$

Theorem 3 implies that, although our test is consistent for any fixed alternative, a cost for this advantage is that it can only detect local alternatives of $O((nm)^{-1/4}) = O(n^{-(1+\beta)/4})$, with $1/2 < \beta < 3/4$ being slightly slower than $n^{-1/2}$. A similar feature is shared with the tests for regression models proposed by Hong and White [17].

REMARK 1. We make several comments on the results.

(i) We make an explicit comparison of (A4), the bandwidth of smoothed periodograms and the rate of convergence with those of related works. By replacing $w^*$ with the integral $m\int_{-1/2}^{1/2} u(x)\,dx$ and rewriting $u(j/m)$ as $u(n\lambda_j/(2\pi m))$, we obtain

$$\frac{w_j}{w^*} \sim \frac{1}{m\int_{-1/2}^{1/2} u(x)\,dx}u\bigg(\frac{n\lambda_j}{2\pi m}\bigg).$$

On the other hand, some authors use $K(\lambda_j/h)/(nh)$ (see, e.g., Paparoditis [29] and [30]) or $M_n K(M_n\lambda_j)/n$ (see, e.g., Eichler [11]) with a kernel function $K(x)$ instead of $w_j/w^*$. $h$ and $M_n$ are called a smoothed bandwidth and an effective number of frequencies, respectively. Then, we have the relation

$$m \sim nh \sim n/M_n$$

as $n \to \infty$, and our $\sqrt{m/n}\,T_n$ corresponds to $T_n$ of [30] and $Q_T$ of [11], respectively. In terms of $m$, (A4) is stronger than the assumptions like $0 < \beta < 1$ (see [29] and [30]) or $1/2 < \beta < 1$ (see [11]), which is the cost for considering nonparametric and semiparametric testing hypotheses in a comprehensive way but not a specific one.

(ii) There are other test statistics alternative to $T_n$. First, in some applications, it may be better to leave out the frequencies where the determinant of $f(\lambda)$ is near zero to make a test statistic more stable. It suggests modifying $T_n$ to

$$T_{n,\phi} = \sum_{t=1}^{[n/2]} \phi_t K(M_t),$$

where $\phi_t = \phi(\lambda_t)$ and $\phi(\lambda)$ is a nonnegative weight function (see also [18]). The limiting behaviors of $T_{n,\phi}$ are obtained in the same way as $T_n$.

Second, we introduce the quadratic function

$$T_{Q,n} = \frac{1}{2} \sum_t \operatorname{tr}[(M_t - I_r)^2]$$

and define the test statistic $\hat{T}_{Q,n}$ by

$$\hat{T}_{Q,n} = \sqrt{\frac{m}{n}} \left( T_{Q,n} - \frac{n}{m}\hat{\eta}_n \right) / \hat{\sigma}_n.$$

A similar idea is proposed by Kakizawa, Shumway and Taniguchi [18] for a discriminant analysis of multivariate time series.

Finally, breaking up the frequency axis with nonoverlapping blocks and using only $M_{(t-1)(m+1)+m/2+1}$ $(t = 1, \ldots, L)$, we define $T_n^*$ by

$$T_n^* = \sum_{t=1}^{L} K(M_{(t-1)(m+1)+m/2+1}),$$

where $L = [n/2]/(m+1)$.

Then, the test statistic is given by

$$\hat{T}_n^* = \frac{m}{L^{1/2}} \left( T_n^* - \frac{2L}{m}\hat{\eta}_n \right) \Big/ \left( \sqrt{B_u/D_u}\,\hat{\sigma}_n \right),$$

where $B_u = (\int_{-1/2}^{1/2} u(x)^2\,dx)^2 / \int_{-1/2}^{1/2} u(x)^4\,dx$. For instance, Wahba [40] applies $T_n^*$ to Example 2 of Section 4 by substituting $K_I(A)$ for $K(A)$.

Similar to the proofs of Theorems 1 and 2, it is shown that the limiting distributions of $\hat{T}_{Q,n}$ and $\hat{T}_n^*$ are $N(0,1)$ under the null hypothesis (1), whereas, if the null hypothesis is false, $\hat{T}_{Q,n}/\sqrt{nm}$ and $\hat{T}_n^*/\sqrt{nm}$ converges in probability to $C_Q^* = \frac{1}{2}\frac{1}{2\pi\sigma^*}\int_0^\pi \operatorname{tr}[(f(\lambda)g(\theta^*, f(\lambda))^{-1} - I_r)^2]\,d\lambda$ and $\sqrt{2D_u/B_u}\,C^*$, respectively.

(iii) For practical use, we need a reasonable criterion to choose a specific test statistic among $K(A)$ mentioned above or their alternatives. Theoretically, various concepts of asymptotic relative efficiency (ARE) of one test relative to the other one have been proposed. They differ from each other in intuitive appeal, or the availability of mathematical tools and efficiency comparison among tests cannot be done in a single ARE (see Serfling [35] for a comprehensive survey of AREs). Hence, the choice of a test statistic among candidates, from a practical point of view, is left to future studies.

Here, we shall give a remark. According to Pitman's ARE approach (see Pitman [32] and Noether [27]), we consider the asymptotic power function when local alternatives converge to the null hypothesis as the sample size goes to infinity. Because $K_I(A), K_J(A)$ and $K_\alpha(A)$ have the same $\kappa_{ab,cd}^{(2)}$ up to the constants, it follows from Theorem 3 that the expectations $\xi/\sigma$ of their limiting distributions under local alternatives are identical to each other, which means that they have the same asymptotic efficiency in Pitman's ARE sense. Similarly, $\hat{T}_{Q,n}$ has the same asymptotic efficiency, whereas $\hat{T}_n^*$ is less efficient than these test statistics, because its expectation of the limiting distribution is $\sqrt{2D_u/B_u}\,\xi/\sigma$ and, by the Schwarz inequality, $2D_u/B_u \le 1$.

**4. Examples.** First, we show some generic formulas, which are helpful for deriving $\eta$ and $\sigma^2$ in Theorem 1. Since

$$(4) \qquad \frac{\partial \breve{g}^{-1}(\lambda)}{\partial y_{\alpha\beta}} = -\breve{g}^{-1}(\lambda)\frac{\partial \breve{g}(\lambda)}{\partial y_{\alpha\beta}}\breve{g}^{-1}(\lambda),$$

we can evaluate $\mu_{\alpha\beta\gamma\nu}(\lambda)$ without calculating $\frac{\partial \breve{g}^{-1}(\lambda)}{\partial y_{\alpha\beta}}$. By applying (4),

$$
\begin{aligned}
&\mu_{\alpha\beta\gamma\nu}(\lambda)\\
&= \frac{1}{2}\operatorname{tr}\left[\frac{\partial \breve{g}(\lambda)}{\partial y_{\alpha\beta}}\breve{g}^{-1}(\lambda)\frac{\partial \breve{g}(\lambda)}{\partial y_{\gamma\nu}}\breve{g}^{-1}(\lambda)\right] - \frac{1}{2}\left[\breve{g}^{-1}(\lambda)\frac{\partial \breve{g}(\lambda)}{\partial y_{\alpha\beta}}\breve{g}^{-1}(\lambda)\right]_{\nu\gamma}\\
&\quad - \frac{1}{2}\left[\breve{g}^{-1}(\lambda)\frac{\partial \breve{g}(\lambda)}{\partial y_{\gamma\nu}}\breve{g}^{-1}(\lambda)\right]_{\beta\alpha} + \frac{1}{2}\breve{g}^{\beta\gamma}(\lambda)\breve{g}^{\nu\alpha}(\lambda)\\
&= \sum_{i=1}^4 \mu_{\alpha\beta\gamma\nu,i}(\lambda), \qquad \text{say.}
\end{aligned}
$$

(5)

Then, for the calculation of $\eta$, we have

$$
\begin{aligned}
&\sum_{\alpha,\beta,\gamma,\nu=1}^r \frac{1}{\pi}\int_0^\pi \mu_{\alpha\beta\gamma\nu,2}(\lambda)\breve{g}_{\alpha\nu}(\lambda)\breve{g}_{\gamma\beta}(\lambda)\,d\lambda\\
&= \sum_{\alpha,\beta,\gamma,\nu=1}^r \frac{1}{\pi}\int_0^\pi \mu_{\alpha\beta\gamma\nu,3}(\lambda)\breve{g}_{\alpha\nu}(\lambda)\breve{g}_{\gamma\beta}(\lambda)\,d\lambda
\end{aligned}
$$

(6)

$$= -\frac{1}{2\pi} \sum_{\alpha,\beta=1}^{r} \int_0^\pi \left[ \frac{\partial \breve{g}(\lambda)}{\partial y_{\alpha\beta}} \right]_{\alpha\beta} d\lambda$$

and

$$(7) \qquad \sum_{\alpha,\beta,\gamma,\nu=1}^{r} \frac{1}{\pi} \int_0^\pi \mu_{\alpha\beta\gamma\nu,4}(\lambda) \breve{g}_{\alpha\nu}(\lambda) \breve{g}_{\gamma\beta}(\lambda) \, d\lambda = \frac{r^2}{2}.$$

Hence, only $\int \mu_{\alpha\beta\gamma\nu,1}(\lambda) \breve{g}_{\alpha\nu}(\lambda) \breve{g}_{\gamma\beta}(\lambda) \, d\lambda$ may requires a laborious evaluation. A similar technique reduces the evaluation of $\sigma^2$ to a simpler one.

Now, we only consider the three examples of those mentioned in Section 1 and proceed to calculate $\eta$ and $\sigma^2$ of them. For simplicity, we assume that $u(x) \equiv 1$ on $[-1/2, 1/2]$ and $w_j \equiv 1$. Then, $C_u = 1/2$ and $D_u = 1/3$.

EXAMPLE 1 (*A separable model*).   The spectral density matrix of a separable stationary process is expressed in the form

$$f(\lambda) = \Sigma \tilde{f}(\lambda),$$

where $\Sigma$ is an $r \times r$ positive definite matrix and $\tilde{f}(\lambda)$ is a scalar-valued nonnegative integrable function in $(-\pi, \pi]$ (see Matsuda and Yajima [25]). Set $v = r^2$ and $\theta = \mathrm{vec}(\Sigma)$, where vec transforms an $r \times r$ matrix into an $r^2$-dimensional vector by stacking the columns of the matrix underneath each other. If we define $g(\theta, y) = \frac{1}{r}(\sum_{\alpha=1}^{r} y_{\alpha\alpha}/\sigma_{\alpha\alpha})\Sigma$, then $f(\lambda) = g(\theta, f(\lambda)) = \frac{1}{r}(\sum_{\alpha=1}^{r} f_{\alpha\alpha}(\lambda)/\sigma_{\alpha\alpha})\Sigma = \Sigma \tilde{f}(\lambda)$.

Next, we observe that $g^{-1}(\theta, y) = r(\sum_{\alpha=1}^{r} y_{\alpha\alpha}/\sigma_{\alpha\alpha})^{-1}\Sigma^{-1}$ and $\partial g(\theta, y)/\partial y_{\alpha\beta} = \delta_{\alpha\beta}/(r\sigma_{\alpha\alpha})\Sigma$, where $\delta_{\alpha\beta}$ is the Kronecker delta. Let $\Sigma_0 = (\sigma_{ab,0})$ and $\Sigma_0^{-1} = (\sigma_0^{ab})$ be the true matrices of $\Sigma$ and $\Sigma^{-1}$, respectively. Then, since $\breve{g}(\lambda) = \Sigma_0 \tilde{f}(\lambda)$, $\partial \breve{g}(\lambda)/\partial y_{\alpha\beta} = \delta_{\alpha\beta}/(r\sigma_{\alpha\alpha,0})\Sigma_0$ and $\breve{g}^{-1}(\lambda) = \Sigma_0^{-1}/\tilde{f}(\lambda)$, from (5),

$$\mu_{\alpha\beta\gamma\nu}(\lambda) = \frac{1}{2\tilde{f}(\lambda)^2} \left( \frac{\delta_{\alpha\beta}\delta_{\gamma\nu}}{r\sigma_{\alpha\alpha,0}\sigma_{\gamma\gamma,0}} - \frac{\delta_{\alpha\beta}\sigma_0^{\nu\gamma}}{r\sigma_{\alpha\alpha,0}} - \frac{\delta_{\gamma\nu}\sigma_0^{\beta\alpha}}{r\sigma_{\gamma\gamma,0}} + \sigma_0^{\beta\gamma}\sigma_0^{\nu\alpha} \right)$$

$$= \sum_{i=1}^{4} \mu_{\alpha\beta\gamma\nu,i}(\lambda), \qquad \text{say.}$$

Thus,

$$\sum_{\alpha\beta\gamma\mu} \frac{1}{\pi} \int_0^\pi \mu_{\alpha\beta\gamma\nu,1}(\lambda) \breve{g}_{\alpha\nu}(\lambda) \breve{g}_{\gamma\beta}(\lambda) \, d\lambda = \frac{\tau}{2r},$$

where $\tau = \sum_{a,b=1}^{r} \sigma_{ab,0}^2/(\sigma_{aa,0}\sigma_{bb,0})$.

Next, from (6),

$$\sum_{\alpha\beta\gamma\mu} \frac{1}{\pi} \int_0^\pi \mu_{\alpha\beta\gamma\nu,2}(\lambda) \breve{g}_{\alpha\nu}(\lambda) \breve{g}_{\gamma\beta}(\lambda) \, d\lambda$$

$$= \sum_{\alpha\beta\gamma\mu} \frac{1}{\pi} \int_0^\pi \mu_{\alpha\beta\gamma\nu,3}(\lambda)\breve{g}_{\alpha\nu}(\lambda)\breve{g}_{\gamma\beta}(\lambda)\,d\lambda$$

$$= -\frac{1}{2r} \sum_{\alpha\beta} \delta_{\alpha\beta}$$

$$= -\frac{1}{2}.$$

Finally, from (7),

$$\eta = \frac{1}{4}\left(\frac{\tau}{r} - 2 + r^2\right).$$

Similarly,

$$\sigma^2 = \frac{1}{6}\left(\frac{\tau^2}{r^2} - 2 + r^2\right).$$

Consistent estimators of $\theta$, $\eta$ and $\sigma^2$ are obtained by substituting $\hat{\Sigma}_n = (\hat{\sigma}_{ab,n}) = \frac{1}{n}\sum_{t=1}^n \mathbf{Z}_t\mathbf{Z}_t'$ for $\Sigma$.

EXAMPLE 2 (*The independence between the component time series*). The independence between the component time series of a stationary Gaussian multivariate time series is equivalent to

$$f_{ab}(\lambda) = 0, \qquad a \neq b,$$

for all $\lambda$ and $a, b = 1, 2, \ldots, r$. Thus, $v = 0$ and $g(y) = \text{diag}(y_{11}, \ldots, y_{rr})$.

Hence, $g^{-1}(y) = \text{diag}(y_{11}^{-1}, \ldots, y_{rr}^{-1})$ and $\partial g(y)/\partial y_{\alpha\beta} = \delta_{\alpha\beta}\,\text{diag}(\delta_{1\alpha}, \ldots, \delta_{r\alpha})$. Then, since $\breve{g}(\lambda) = \text{diag}(f_{11}(\lambda), \ldots, f_{rr}(\lambda))$, $\partial\breve{g}(\lambda)/\partial y_{\alpha\beta} = \delta_{\alpha\beta}\,\text{diag}(\delta_{1\alpha}, \ldots, \delta_{r\alpha})$ and $\breve{g}^{-1}(\lambda) = \text{diag}(f_{11}^{-1}(\lambda), \ldots, f_{rr}^{-1}(\lambda))$, from (5),

$$\mu_{\alpha\beta\gamma\nu}(\lambda) = \frac{1}{2}\left(\frac{\delta_{\alpha\nu}\delta_{\beta\gamma}}{f_{\alpha\alpha}(\lambda)f_{\beta\beta}(\lambda)} - \frac{\delta_{\alpha\beta}\delta_{\alpha\gamma}\delta_{\gamma\nu}}{f_{\alpha\alpha}(\lambda)^2}\right).$$

Thus, $\eta$ and $\sigma^2$ are the known constants expressed as

$$\eta = \tfrac{1}{4}(r^2 - r), \qquad \sigma^2 = \tfrac{1}{6}(r^2 - r).$$

Hong [16] and Eichler [11] consider the same problem for the bivariate time series including non-Gaussian time series. Our result is a generalization to a multivariate one, though being derived under Gaussian assumption.

Finally, we consider an example where a nonlinear constraint is imposed on the null hypothesis.

EXAMPLE 3 (*The conditional independence*). Set $Z_a = \{Z_{at}, -\infty < t < \infty\}, Y_{ab} = \{Y_{ab,t}, -\infty < t < \infty\}$ where $Y_{ab,t} = \{Z_{jt}, j \neq a, b\}$ is an $(r-2)$-dimensional random vector. Then, the conditional independence between $Z_a$ and $Z_b$ given $Y_{ab}$ is defined by

$$(8) \qquad \mathrm{Cov}(\varepsilon_{a|\{a,b\}^c}(s), \varepsilon_{b|\{a,b\}^c}(t)) = 0$$

for all $s, t \in \mathcal{Z} = \{0, \pm 1, \pm 2, \ldots\}$, where

$$\varepsilon_{a|\{a,b\}^c}(t) = Z_{at} - \sum_{u=-\infty}^{\infty} d_a^*(t-u)' Y_{ab,u},$$

and $d_a^*(t-u)$ is the $(r-2)$-dimensional vector which minimizes

$$E\left(Z_{at} - \sum_{u=-\infty}^{\infty} d_a(t-u)' Y_{ab,u}\right)^2.$$

Then, the relation (8) is equivalent to

$$f^{ab}(\lambda) = 0, \qquad -\pi \leq \lambda \leq \pi$$

(see, e.g., Dahlhaus [7]).

Now, let $V = \{1, \ldots, r\}$ and $E$ be a subset of $V \times V$. Consider the null hypothesis

$$f^{ab}(\lambda) = 0, \qquad (a, b) \notin E$$

for all $\lambda \in (-\pi, \pi]$. Then, $v = 0$ and

$$(9) \qquad \begin{aligned} g_{ab}(y) &= y_{ab}, & (a, b) \in E, \\ g^{ab}(y) &= 0, & (a, b) \notin E. \end{aligned}$$

A graph $G = (V, E)$ is the partial correlation graph of a time series $\mathbf{Z}_t$ if $(a, b) \notin E \Leftrightarrow f^{ab}(\lambda) = 0, \forall \lambda \in (-\pi, \pi]$ (see Dahlhaus [7]). It is proved by the implicit function theorem, in the same way as Lemma 7 of Matsuda, Yajima and Tong [26], that the constraint (9) is expressed in the form (1).

Next, we show the outline of the derivations of $\eta$ and $\sigma^2$, because it is given in the same way as Theorem 3 of Matsuda, Yajima and Tong [26] and Theorem 2 of Matsuda [24] by applying Lemma 8 of Matsuda, Yajima and Tong [26].

Let $M = \#\{(a, b) | (a, b) \notin E, a < b\}$ and $I_C$ be the indicator function which is 1 if the condition $C$ is true and 0 otherwise. Then, we have

$$\begin{aligned}
\mu_{\alpha\beta\gamma\nu,1}&(\lambda) \\
&= \frac{1}{2} I_{(\alpha,\beta) \in E} I_{(\gamma,\nu) \in E} \sum_{ab} \left( \sum_e I_{(a,e) \notin E} \left[ \frac{\partial \breve{g}(\lambda)}{\partial y_{\alpha\beta}} \right]_{ae} \breve{g}^{eb}(\lambda) + I_{a=\alpha} \breve{g}^{\beta b}(\lambda) \right) \\
&\qquad\qquad \times \left( \sum_{e'} I_{(b,e') \notin E} \left[ \frac{\partial \breve{g}(\lambda)}{\partial y_{\gamma\nu}} \right]_{be'} \breve{g}^{e'a}(\lambda) + I_{b=\gamma} \breve{g}^{\nu a}(\lambda) \right).
\end{aligned}$$

Thus,

$$\sum_{\alpha,\beta,\gamma,\nu=1}^{r} \frac{1}{\pi} \int_0^{\pi} \mu_{\alpha\beta\gamma\nu,1}(\lambda) \breve{g}_{\alpha\nu}(\lambda) \breve{g}_{\gamma\beta}(\lambda) \, d\lambda$$

$$= -M + \frac{1}{2\pi} \int_0^{\pi} [\mathrm{tr}(I_r)]^2 \, d\lambda$$

$$= -M + \frac{1}{2} r^2.$$

From (6),

$$\sum_{\alpha,\beta,\gamma,\nu=1}^{r} \frac{1}{\pi} \int_0^{\pi} \mu_{\alpha\beta\gamma\nu,2}(\lambda) \breve{g}_{\alpha\nu}(\lambda) \breve{g}_{\gamma\beta}(\lambda) \, d\lambda = -\frac{1}{2}(r^2 - 2M).$$

Then, it follows from (7) that $\eta = M/2$. Similarly, $\sigma^2 = M/3$.

REMARK 2. All of the $g(\theta, y)$ in Examples 1–3 map the space of positive definite matrices into itself (see Dempster [9] for Example 3). However, it is not necessary for $g(\theta, y)$ to satisfy this condition under alternatives, because Theorems 1 and 3 depend only on the behavior of $g(\theta, f(\lambda))$ in the neighborhood of the null hypothesis, and Theorem 2 still holds by assigning a sufficiently large value to $\hat{T}_n$ if $\hat{f}_{R,t}$ is not positive definite, and, consequently, $K(M_t)$ cannot be defined. For practical use, being not positive definite, $\hat{f}_{R,t}$ gives strong evidence against the null hypothesis, and, consequently, we can reject it without causing any serious problem.

**5. Simulation results.** We conduct some computational simulations to see the performance of the test statistics. Consider the following three-dimensional model for testing the independence between the component series:

$$(10) \qquad \mathbf{Z}_t = \begin{pmatrix} 0.7 & \phi & 0.0 \\ 0.0 & -0.5 & \phi \\ 0.0 & 0.0 & 0.6 \end{pmatrix} \mathbf{Z}_{t-1} + \varepsilon_t,$$

where $\varepsilon_t$ are independent normal variables with mean 0 and covariance matrix $I_3$. The component series of (10) are mutually independent if $\phi = 0.0$ and are dependent otherwise.

Applying our test statistic $\hat{T}_n$ with $K_I(A)$ and $u(x) \equiv 1$ on $[-1/2, 1/2]$, we test the independence between the component series mentioned in Example 2 of Section 4. We also considered Wahba's test statistic $\hat{T}_n^*$ and the quadratic test statistic $\hat{T}_{Q,n}$ for comparison. We examine performances of the tests under the null and alternative hypotheses in Tables 1 and 2, respectively.

Table 1 shows mean, variance, skewness, kurtosis, 5% upper quantile of the null distributions and the empirical sizes under 5% asymptotic significance level based on the 1000 replications of the process (10) with $\phi = 0$. We find that all the quantities except for skewness converge to the asymptotic limits relatively fast. However, it should be pointed out that bias exists in every test statistic in the small sample size, which causes nonnegligible size distortions. Hence, for fair power comparison, all the sizes of the test statistics are adjusted to be exactly 0.05 in the same way as Haug [15] and Saikkonen and Luukkonen [34]. Table 2 shows the empirical powers of the three statistics based on the 1000 replications of the process (10) with $\phi = 0.1$ and $0.2$. We show the results when the bandwidth is predetermined as $n = 101, m = 16$ and $n = 201, m = 30$ and is selected by cross validation,

TABLE 1
*Comparison of the null distribution of the tests for independence where 5% size is the empirical frequency of the rejection by the asymptotic 5% critical point*

| Test | Sample size | $m$ | Mean | Var | Skew | Kurt | 5% qtl | 5% size |
|------|------------|-----|------|-----|------|------|--------|---------|
| $\hat{T}_n$ | 101 | 16 | 0.163 | 1.04 | 0.66 | 4.01 | 1.94 | 0.075 |
| | 201 | 30 | 0.097 | 1.02 | 0.62 | 3.61 | 1.94 | 0.083 |
| | 501 | 60 | 0.073 | 0.94 | 0.59 | 3.62 | 1.74 | 0.059 |
| | 1001 | 120 | 0.069 | 0.94 | 0.44 | 3.09 | 1.75 | 0.064 |
| $\hat{T}_{Q,n}$ | 101 | 16 | $-0.163$ | 0.78 | 0.59 | 3.70 | 1.44 | 0.036 |
| | 201 | 30 | $-0.084$ | 0.87 | 0.59 | 3.56 | 1.63 | 0.048 |
| | 501 | 60 | $-0.028$ | 0.87 | 0.59 | 3.63 | 1.60 | 0.045 |
| | 1001 | 120 | 0.020 | 0.90 | 0.42 | 3.06 | 1.67 | 0.053 |
| $\hat{T}_n^*$ | 101 | 16 | 0.161 | 1.14 | 0.76 | 3.91 | 2.07 | 0.095 |
| | 201 | 30 | 0.073 | 1.00 | 0.69 | 3.43 | 1.91 | 0.075 |
| | 501 | 60 | 0.067 | 1.01 | 0.64 | 3.41 | 1.80 | 0.070 |
| | 1001 | 120 | 0.050 | 1.01 | 0.47 | 3.08 | 1.89 | 0.072 |
| | $\infty$ | | 0.0 | 1.0 | 0.0 | 3.0 | 1.64 | 0.050 |

TABLE 2
*Power comparisons of the size adjusted tests for independence under 5% empirical significance level*

| Sample size | $\phi$ | $\hat{T}_n$ | $\hat{T}_{Q,n}$ | $\hat{T}_n^*$ | $\hat{T}_n$ | $\hat{T}_{Q,n}$ | $\hat{T}_n^*$ |
|-------------|--------|-------------|-----------------|---------------|-------------|-----------------|---------------|
| | | | $m = 16$ | | | CVLL | |
| 101 | 0.1 | 0.130 | 0.118 | 0.097 | 0.110 | 0.118 | 0.098 |
| | 0.2 | 0.503 | 0.482 | 0.275 | 0.451 | 0.451 | 0.276 |
| | | | $m = 30$ | | | CVLL | |
| 201 | 0.1 | 0.231 | 0.231 | 0.163 | 0.206 | 0.232 | 0.145 |
| | 0.2 | 0.858 | 0.859 | 0.719 | 0.831 | 0.840 | 0.621 |

which minimizes

$$\text{CVLL}(m) = \frac{1}{n} \sum_{j=1}^{[n/2]} (\text{tr}(I_{Z,j}\hat{f}_{U,j,-j}^{-1}) + \log \det(\hat{f}_{U,j,-j})),$$

where

$$\hat{f}_{U,j,-j} = \frac{1}{m} \sum_{k=-m/2, k\neq 0}^{m/2} I_{Z,j+k}.$$

Alternatives to CVLL can be the methods proposed by Lee [21] and Ombao et al. [28]. For tentative comparison between the performance of the fixed bandwidth selection and the data driven one, we adopt CVLL. However, we should remark that CVLL can select the bandwidth minimizing the mean squared error asymptotically if we put $\beta = 4/5$, which lies outside the interval of (A4) (see, e.g., Beltrao and Bloomfield [1], Matsuda and Yajima [25]). This issue often emerges in nonparametric hypothesis testing (see, e.g., Zhang [41] and the references therein and Fan and Yao [12], Section 9.2.7) and more rigorous consideration is left to future studies.

We find from Table 2 that our statistic has almost the same power as the quadratic one, whereas the Wahba's is significantly less powerful than the others, which reinforces that Wahba's statistic is less efficient than the others in Pitman's ARE sense as mentioned in Remark 1(iii), because $2D_u/C_u = 2/3$ is less than 1.

## 6. Proofs of theorems.

PROOF OF THEOREM 1. Applying the Taylor expansion and Lemma 2 to $K(M_t)$ and noting that the first-order terms of the expansion vanish, we obtain

$$
\begin{aligned}
T_n = {}& \frac{1}{2} \sum_{i=1}^{[n/2]} \text{tr}[(M_t - I_r)^2] \\
& + \frac{1}{6} \sum_{a,b,c,d,e,f=1}^{r} k_{ab,cd,ef}^{(3)} \\
& \qquad\qquad \times \sum_{t=1}^{[n/2]} (m_{ab,t} - \delta_{ab})(m_{cd,t} - \delta_{cd})(m_{ef,t} - \delta_{ef}) \\
& + o_p(1).
\end{aligned}
\tag{11}
$$

It is shown later that the second term on the right-hand side of (11) is $o_p(1)$. Hence, we shall consider the limiting distribution of the first term.

Note that $f(\lambda) = \breve{g}(\lambda)$ if the null hypothesis (1) is true. Then, it follows from Lemma 2 and $\hat{\theta}_n - \theta_0 = O_p(n^{-1/2})$, by the Taylor expansion, that

$$
\begin{aligned}
M_t &- I_r \\
&= \hat{f}_{U,t}(\hat{f}_{R,t}^{-1} - \breve{g}_t^{-1}) + (\hat{f}_{U,t} - \breve{g}_t)\breve{g}_t^{-1} \\
&= \breve{g}_t\{g^{-1}(\hat{\theta}_n, \hat{f}_{U,t}) - g^{-1}(\theta_0, f_t)\} \\
&\quad + (\hat{f}_{U,t} - \breve{g}_t)\{g^{-1}(\hat{\theta}_n, \hat{f}_{U,t}) - g^{-1}(\theta_0, f_t)\} + (\hat{f}_{U,t} - \breve{g}_t)\breve{g}_t^{-1} \\
(12) \quad &= \sum_{\alpha,\beta=1}^{r} \breve{g}_t(\partial \breve{g}_t^{-1}/\partial y_{\alpha\beta})(\hat{f}_{U,\alpha\beta,t} - f_{\alpha\beta,t}) + \sum_{i=1}^{v} \breve{g}_t(\partial \breve{g}_t^{-1}/\partial\theta_i)(\hat{\theta}_{in} - \theta_{i0}) \\
&\quad + \frac{1}{2}\sum_{\alpha,\beta,\gamma,\nu=1}^{r} \breve{g}_t(\partial^2 \breve{g}_t^{-1}/\partial y_{\alpha\beta}\partial y_{\gamma\nu})(\hat{f}_{U,\alpha\beta,t} - f_{\alpha\beta,t})(\hat{f}_{U,\gamma\nu,t} - f_{\gamma\nu,t}) \\
&\quad + o_p(n^{-3/4}) + (\hat{f}_{U,t} - \breve{g}_t)\{g^{-1}(\hat{\theta}_n, \hat{f}_{U,t}) - g^{-1}(\theta_0, f_t)\} \\
&\quad + (\hat{f}_{U,t} - \breve{g}_t)\breve{g}_t^{-1}.
\end{aligned}
$$

By applying Lemmas 2, 6 and 7 and $\hat{\theta}_n - \theta_0 = O_p(n^{-1/2})$ to (12),

$$
\begin{aligned}
\sum_t &\operatorname{tr}[(M_t - I_r)^2] \\
(13) \quad &= \sum_t \operatorname{tr}\left[\left(\sum_{\alpha\beta} \breve{g}_t(\partial \breve{g}_t^{-1}/\partial y_{\alpha\beta})(\hat{f}_{U,\alpha\beta,t} - f_{\alpha\beta,t}) + (\hat{f}_{U,t} - \breve{g}_t)\breve{g}_t^{-1}\right)^2\right] \\
&\quad + o_p(n^{1/2}/m^{1/2}) \\
&= T_n^{(1)} + o_p(n^{1/2}/m^{1/2}) \qquad \text{say.}
\end{aligned}
$$

Next, define $T_n^{(2)}$ by replacing $\hat{f}_{U,t} - \breve{g}_t$ and $\hat{f}_{U,t} - f_t$ of $T_n^{(1)}$ with $\hat{f}_t^{\varepsilon} - E(\hat{f}_t^{\varepsilon})$ of Section 7. Then, from Lemma 8, it suffices to derive the limiting distribution of $\frac{1}{2}T_n^{(2)}$.

By an elementary calculation,

$$
\frac{1}{2}T_n^{(2)} = \sum_{\alpha,\beta,\gamma,\nu=1}^{r} \sum_t \mu_{\alpha\beta\gamma\nu,t,n}\hat{y}_{\alpha\beta,t}^{\varepsilon}\hat{y}_{\gamma\nu,t}^{\varepsilon},
$$

where $\mu_{\alpha\beta\gamma\nu,t,n} = \mu_{\alpha\beta\gamma\nu}(\lambda_t)$. Then, if we define $\tilde{T}_n$ by

$$
\tilde{T}_n = \sqrt{\frac{m}{n}}\left(\frac{1}{2}T_n^{(2)} - \frac{n}{m}\eta\right),
$$

it follows from Lemmas 9, 10 and 11 that $\lim_n E(\tilde{T}_n) = 0$, $\lim_n \operatorname{Var}(\tilde{T}_n) = \sigma^2$ and for $k \geq 3$, the $k$th order cumulant of $\tilde{T}_n$ converges to 0 as $n \to \infty$. Hence, the limiting distribution of $\tilde{T}_n$ is $N(0, \sigma^2)$.

Finally, by the same argument as (12), it is shown that the second term on the right-hand side of (11) is negligible.  $\square$

PROOF OF THEOREM 2.    $\hat{T}_n$ is written as

$$\sqrt{nm}\left(\frac{T_n}{n} - \frac{1}{m}\hat{\eta}_n\right)/\hat{\sigma}_n.$$

Then, it follows from Lemma 2 and Exercise 1.7.4 of Brillinger [2] that $T_n/(n\hat{\sigma}_n)$ converges in probability to $\frac{1}{2\pi\sigma^*}\int_0^\pi K(f(\lambda)g(\theta^*, f(\lambda)^{-1})\,d\lambda$ as $n \to \infty$. Thus, the assertion is shown immediately.  $\square$

PROOF OF THEOREM 3.    First, we note that, under $H_{an}$, Lemmas 1–11 of Section 7 are still true if we substitute $f_{t,n}$ for $f_t$ in the definitions of $y_t$ and $y_t^\varepsilon$ and add the term of $O(1/(mn)^{1/4})$ to the right-hand side of Lemma 9. Next, since $M_t - I_r = o_p(n^{-1/4})$ is also true under $H_{an}$, analogous to (13), we have

$$\sum_t \text{tr}[(M_t - I_r)^2]$$

$$= \sum_t \text{tr}\left[\left(\sum_{\alpha\beta} \check{g}_{t,n}(\partial \check{g}_{t,n}^{-1}/\partial y_{\alpha\beta})(\hat{f}_{U,\alpha\beta,t} - f_{\alpha\beta,t,n}) + (\hat{f}_{U,t} - f_{t,n})\check{g}_{t,n}^{-1}\right)^2\right]$$

$$+ \sum_t \text{tr}\left[\left((f_t - \check{g}_{t,n})\check{g}_{t,n}^{-1} + \frac{1}{(nm)^{1/4}}f_t^*\check{g}_{t,n}^{-1}\right)^2\right] + o_p(n^{1/2}/m^{1/2})$$

$$= \overline{T}_n^{(1)} + \xi_n + o_p(n^{1/2}/m^{1/2}),\qquad \text{say.}$$

Then, by the Taylor expansion and Exercise 1.7.14 of Brillinger [2],

$$\xi_n = \frac{n^{1/2}}{m^{1/2}}(2\xi) + o(n^{1/2}/m^{1/2}).$$

Hence,

$$\hat{T}_n = \sqrt{\frac{m}{n}}\left(T_n - \frac{n}{m}\hat{\eta}_n\right)/\hat{\sigma}_n$$

(14)
$$= \sqrt{\frac{m}{n}}\left(\frac{1}{2}\overline{T}_n^{(1)} - \frac{n}{m}\eta\right)/\hat{\sigma}_n$$

$$+ \xi/\hat{\sigma}_n + \sqrt{\frac{n}{m}}\{(\eta - \eta_n) + (\eta_n - \hat{\eta}_n)\}/\hat{\sigma}_n + o_p(1).$$

The limiting distribution of the first term on the right-hand side of (14) is the standard normal distribution. Then, the assertion is obtained, because $\eta_n - \eta = O(1/(mn)^{1/4})$ and $\sigma_n^2 - \sigma^2 = O(1/(mn)^{1/4})$.  $\square$

**7. Lemmas.** First, we introduce some random variables, notation and technical remarks. Define $y_t = (y_{ab,t}), \hat{y}_t = (\hat{y}_{ab,t}), t = 1, \ldots, [n/2]$, by $y_t = \hat{f}_{U,t} - f_t$ and $\hat{y}_t = \hat{f}_{U,t} - E(\hat{f}_{U,t})$, respectively.

(A2) and (A3) assure that $\int_{-\pi}^{\pi} \log \det f(\lambda) \, d\lambda > -\infty$. Then, under (A1), $\mathbf{Z}_t$ is expressed as

$$\mathbf{Z}_t = \sum_{j=0}^{\infty} \Phi_j \varepsilon_{t-j},$$

where $\Phi_j$ is an $r \times r$ matrix and $\sum_{j=0}^{\infty} \operatorname{tr}(\Phi_j \Phi_j') < \infty$ and $\varepsilon_t = (\varepsilon_{1t}, \ldots, \varepsilon_{rt})'$ is a mutually independent zero-mean Gaussian process with covariance matrix $I_r$. Then, denote the discrete Fourier transform of $\varepsilon_{at}$ by $W_a^\varepsilon(\lambda) = \frac{1}{\sqrt{2\pi n}} \sum_{t=1}^n \varepsilon_{at} \times \exp(it\lambda)$ for $a = 1, \ldots, r$ and the cross periodogram of $\varepsilon_{at}$ and $\varepsilon_{bt}$ and the periodogram matrix by $I_{ab}^\varepsilon(\lambda) = W_a^\varepsilon(\lambda)\overline{W_b^\varepsilon(\lambda)}$ and $I^\varepsilon(\lambda) = (I_{ab}^\varepsilon(\lambda))$, respectively.

Next, define the $r \times r$ matrix $\hat{f}_t^\varepsilon = (\hat{f}_{ab,t}^\varepsilon)$ by

$$\hat{f}_{ab,t}^\varepsilon = \frac{1}{w^*} \sum_{j=-m/2}^{m/2} w_j \Phi_a(\exp(i\lambda_{t+j})) I_{t+j}^\varepsilon \Phi_b'(\exp(-i\lambda_{t+j})), \qquad 1 \le t \le [n/2],$$

where $I_j^\varepsilon = I^\varepsilon(\lambda_j)$ and $\Phi_a(e^{i\lambda})$ is the $a$th row vector of $\Phi(e^{i\lambda}) = \sum_{j=0}^{\infty} \Phi_j e^{ij\lambda}$.

Then, define $y_t^\varepsilon = \hat{f}_t^\varepsilon - f_t$ and $\hat{y}_t^\varepsilon = \hat{f}_t^\varepsilon - E(\hat{f}_t^\varepsilon)$, respectively. We shall derive the main results by showing that the limiting behaviors of $\hat{T}_n$ remain unchanged if we substitute $y_t^\varepsilon(\hat{y}_t^\varepsilon)$ for $y_t(\hat{y}_t)$. $y_t^\varepsilon$ and $\hat{y}_t^\varepsilon$ are more tractable than $y_t$ and $\hat{y}_t$, because $2\pi I_j^\varepsilon, j = 0, 1, \ldots, [n/2]$, are mutually independent random variables with Wishart distributions and, consequently, $y_t^\varepsilon(\hat{y}_t^\varepsilon)$ and $y_s^\varepsilon(\hat{y}_s^\varepsilon)$ are mutually independent for $|t - s| > m$.

LEMMA 1. *Under* (A3),

$$E(W_{a,j}\overline{W}_{b,j}) - f_{ab,j} = O\left(\frac{\log n}{n}\right), \qquad -m/2 + 1 \le j \le [n/2] + m/2,$$

$$E(W_{a,j}W_{b,k}) = O\left(\frac{\log n}{n}\right), \qquad -m/2 + 1 \le j \le k \le [n/2] + m/2,$$
$$j + k \ne 0, n,$$

$$E(W_{a,j}W_{b,k}) - f_{ab,j} = O\left(\frac{\log n}{n}\right), \qquad -m/2 + 1 \le j \le k \le [n/2] + m/2,$$
$$j + k = 0, n,$$

$$E(W_{a,j}\overline{W}_{b,k}) = O\left(\frac{\log n}{n}\right), \qquad -m/2 + 1 \le j < k \le [n/2] + m/2,$$

*uniformly in $j$ and $k$ for $a, b = 1, 2, \ldots, r$.*

Proof.   The assertion is shown in the same way as Lemma 1 of Matsuda and Yajima [25] by noting that $W_{a,j}W_{b,k} = W_{a,j}\overline{W}_{b,j}$ for $j + k = 0, n$.   □

Lemma 2.   *Under* (A1)–(A4),

$$\sup_{t=\pm 1,\ldots,\pm[n/2]} |\hat{f}_{U,ab,t} - f_{ab,t}| = o_p(n^{-1/4}).$$

Proof.   Applying Lemma 1, we see that Lemmas 2 and 3 of Matsuda and Yajima [25] are still true for $\hat{f}_{U,ab,t}$ and $\hat{f}_{ab,t}^{\varepsilon}$. Then, the assertion is proved in the same way as Proposition 1 of Matsuda and Yajima [25].   □

Lemma 3.   *Under* (A3),

$$E(y_{ab,t}) = O(m^2/n^2),$$

*uniformly in* $t = 1, 2, \ldots, [n/2]$, *for* $a, b = 1, 2, \ldots, r$.

Proof.   It follows from Lemma 1, by the Taylor expansion, that

$$E(y_{ab,t}) = \frac{1}{w^*} \sum_{j=-m/2}^{m/2} w_j(f_{ab,t+j} + O(\log n/n)) - f_{ab,t}$$

$$= \frac{1}{w^*} \sum_{j=-m/2}^{m/2} w_j\left(f_{ab,t} + \frac{2\pi j}{n}f'_{ab,t} + O(j^2/n^2)\right) + O(\log n/n) - f_{ab,t}$$

$$= O(m^2/n^2) + O(\log n/n)$$

$$= O(m^2/n^2),$$

uniformly in $t$ where $f'_{ab,t} = df_{ab}(\lambda)/d\lambda|_{\lambda=\lambda_t}$.   □

Lemma 4.   *Under* (A1), (A3) *and* (A4),

$$\mathrm{Cov}(y_{ab,t}, y_{cd,s}) = \begin{cases} O(1/m), & \text{if } |t-s| \leq m, \\ O(\log^2 n/n^2), & \text{if } |t-s| > m, \end{cases}$$

*uniformly in* $t, s = 1, 2, \ldots, [n/2]$, *for* $a, b, c, d = 1, 2, \ldots, r$.

Proof.   We observe that

$$\mathrm{Cov}(y_{ab,t}, y_{cd,s})$$

$$= E(y_{ab,t}\overline{y}_{cd,s}) - E(y_{ab,t})E(\overline{y}_{cd,s})$$

$$= \frac{1}{(w^*)^2} \sum_{j,k=-m/2}^{m/2} w_j w_k[E(W_{a,t+j}\overline{W}_{c,s+k})E(\overline{W}_{b,t+j}W_{d,s+k})$$

$$+ E(W_{a,t+j}W_{d,s+k})E(\overline{W}_{b,t+j}\overline{W}_{c,s+k})].$$

Then, the assertion follows immediately from Lemma 1 by noting that $t+j \neq s+k$ and $t+j+s+k \neq 0, n$ for any $j$ and $k$ if $|t-s| > m$. □

LEMMA 5. *Under* (A1), (A3) *and* (A4),

$$\text{cum}(y_{a_1 b_1, t_1}, y_{a_2 b_2, t_2}, \ldots, y_{a_k b_k, t_k}) = O(m^{1-k}),$$

*uniformly in* $1 \leq t_i \leq [n/2]$ *for any* $k \geq 3$ *and* $1 \leq a_i, b_i \leq r$, $i = 1, \ldots, k$.

PROOF. The assertion is obtained, in the same way as the proof of Theorem 7.4.4 of Brillinger [2], by applying Lemma 1 and Theorems 2.3.1 and 2.3.2 of [2]. □

From now on, let $\mu_{t,n}, t = 1, \ldots, [n/2], n = 1, 2, \ldots$, be constants bounded in $t$ and $n$.

LEMMA 6. *Under* (A1), (A3) *and* (A4),

$$\sum_{t=1}^{[n/2]} \mu_{t,n} y_{ab,t} = O_p(n^{1/2})$$

*for* $a, b = 1, \ldots, r$.

PROOF. We may assume, without loss of generality, that $\mu_{t,n} \equiv 1$. Then, it follows from Lemmas 3 and 4 that

$$E \left| \sum_t y_{ab,t} \right|^2 = \sum_{s,t} \text{Cov}(y_{ab,t}, y_{ab,s}) + \left| \sum_t E(y_{ab,t}) \right|^2$$

$$= O(n) + O(\log^2 n) + O(m^4/n^2)$$

$$= O(n),$$

which implies the assertion. □

LEMMA 7. *Under* (A1), (A3) *and* (A4),

$$\sum_{t=1}^{[n/2]} \mu_{t,n} y_{ab,t} y_{cd,t} y_{ef,t} = o_p(1)$$

*for* $a, b, c, d, e, f = 1, \ldots, r$.

PROOF. We may assume, without loss of generality, that $\mu_{t,n} \equiv 1$. It follows from Lemmas 2 and 3 that

$$\hat{y}_{ab,t} = y_{ab,t} + O(m^2/n^2) = o_p(n^{-1/4}),$$

uniformly in $t$. Thus,

$$(15) \qquad \sum_t y_{ab,t} y_{cd,t} y_{ef,t} = \sum_t \hat{y}_{ab,t} \hat{y}_{cd,t} \hat{y}_{ef,t} + o_p(1).$$

Now, we evaluate the first term on the right-hand side of (15). We have

$$(16) \qquad E \left| \sum_t \hat{y}_{ab,t} \hat{y}_{cd,t} \hat{y}_{ef,t} \right|^2$$

$$= \sum_{s,t} \text{Cov}(\hat{y}_{ab,t} \hat{y}_{cd,t} \hat{y}_{ef,t}, \hat{y}_{ab,s} \hat{y}_{cd,s} \hat{y}_{ef,s}) + \left| \sum_t E(\hat{y}_{ab,t} \hat{y}_{cd,t} \hat{y}_{ef,t}) \right|^2.$$

Noting $E(\hat{y}_{ab,t}) = 0$ and applying Theorems 2.3.1 and 2.3.2 of Brillinger [2], we observe that the first term on the right-hand side of (16) is equal to

$$\sum_{s,t} \text{cum}(\hat{y}_{ab,t} \hat{y}_{cd,t} \hat{y}_{ef,t}, \hat{y}_{ba,s} \hat{y}_{dc,s} \hat{y}_{fe,s})$$

$$(17) \qquad \begin{aligned} = \sum_{s,t} [ &\text{cum}(y_{ab,t}, y_{cd,t}, y_{ef,t}, y_{ba,s}, y_{dc,s}, y_{fe,s}) \\ &+ \text{cum}(y_{ab,t}, y_{cd,t}) \text{cum}(y_{ef,t}, y_{ba,s}, y_{dc,s}, y_{fe,s}) \\ &+ \text{cum}(y_{ab,t}, y_{ba,s}) \text{cum}(y_{cd,t}, y_{ef,t}, y_{dc,s}, y_{fe,s}) \\ &+ \text{cum}(y_{ab,t}, y_{cd,t}, y_{ba,s}) \text{cum}(y_{ef,t}, y_{dc,s}, y_{fe,s}) \\ &+ \text{cum}(y_{ab,t}, y_{cd,t}) \text{cum}(y_{ba,s}, y_{dc,s}) \text{cum}(y_{ef,t}, y_{fe,s}) \\ &+ \text{cum}(y_{ab,t}, y_{ba,s}) \text{cum}(y_{cd,t}, y_{dc,s}) \text{cum}(y_{ef,t}, y_{fe,s}) \\ &+ \text{the remainder terms}]. \end{aligned}$$

Any term of the remainder ones on the right-hand side of (17) is expressed in the same form as one of the preceding six terms. Hence, it suffices to evaluate these terms. It follows from Lemmas 4 and 5 that the first one is $O(n^2/m^5)$, the second and fourth ones are $O(n^2/m^4)$, the third one is $O(n/m^3)$ and the fifth and sixth ones are $O(n/m^2)$, respectively.

Similarly, the second term is equal to

$$\left| \sum_t \text{cum}(\hat{y}_{ab,t}, \hat{y}_{cd,t}, \hat{y}_{ef,t}) \right|^2 = \left| \sum_t \text{cum}(y_{ab,t}, y_{cd,t}, y_{ef,t}) \right|^2$$

$$= O(n^2/m^4).$$

Thus, the proof is completed.   $\square$

LEMMA 8.   *Under* (A1)–(A4),

$$T_n^{(1)} - T_n^{(2)} = o_p(n^{1/2}/m^{1/2}).$$

PROOF. It suffices to show that

$$\sum_{t=1}^{[n/2]} \mu_{t,n}(y_{\alpha\beta,t} y_{\gamma\nu,t} - \hat{y}_{\alpha\beta,t}^{\varepsilon} \hat{y}_{\gamma\nu,t}^{\varepsilon}) = o_p(n^{1/2}/m^{1/2}) \tag{18}$$

for $\alpha, \beta, \gamma, \nu = 1, \ldots, r$, because $T_n^{(1)} - T_n^{(2)}$ is a summation of the terms expressed in the same form as (18). We may assume that $\mu_{t,n} \equiv 1$ without loss of generality. By the Schwarz inequality, the absolute value of the term on the left-hand side of (18) is bounded by

$$\left(\sum_t |y_{\alpha\beta,t} - \hat{y}_{\alpha\beta,t}^{\varepsilon}|^2\right)^{1/2} \left(\sum_t |y_{\gamma\nu,t}|^2\right)^{1/2}$$
$$+ \left(\sum_t |\hat{y}_{\alpha\beta,t}^{\varepsilon}|^2\right)^{1/2} \left(\sum_t |y_{\gamma\nu,t} - \hat{y}_{\gamma\nu,t}^{\varepsilon}|^2\right)^{1/2}. \tag{19}$$

The first term of (19) is bounded by

$$\left(2\sum_t |y_{\alpha\beta,t} - y_{\alpha\beta,t}^{\varepsilon}|^2 + 2\sum_t |y_{\alpha\beta,t}^{\varepsilon} - \hat{y}_{\alpha\beta,t}^{\varepsilon}|^2\right)^{1/2} \left(\sum_t |y_{\gamma\nu,t}|^2\right)^{1/2}. \tag{20}$$

It follows, from Lemmas 3 and 4 and the assertion during the proof of Lemma 2 of Matsuda and Yajima [25], that the first term of (20) is equal to $O_p((\log n/m)^{1/2}) + O(m^2/n^{3/2})$ and the second term is $O_p((n/m)^{1/2})$. The second term of (19) is evaluated in the same way. Thus, the proof is complete. □

LEMMA 9. *Let $\mu(\lambda)$ be a differentiable function in $[0, \pi]$, and set $\mu_{t,n} = \mu(\lambda_t), t = 1, \ldots, [n/2]$.*

*Then, under* (A1)–(A4),

$$\frac{m}{n} \sum_{t=1}^{[n/2]} E(\mu_{t,n} \hat{y}_{ab,t}^{\varepsilon} \hat{y}_{cd,t}^{\varepsilon}) = \frac{C_u}{\pi} \int_0^{\pi} \mu(\lambda) f_{ad}(\lambda) f_{cb}(\lambda) \, d\lambda + O(m/n)$$

*for $a, b, c, d = 1, \ldots, r$ where*

$$C_u = \frac{1}{2} \int_{-1/2}^{1/2} u(x)^2 \, dx \bigg/ \left(\int_{-1/2}^{1/2} u(x) \, dx\right)^2.$$

PROOF. Applying Lemma 1 to $W_{a,j}^{\varepsilon}, a = 1, \ldots, r$, and noting that the right-hand side terms of the lemma are exactly equal to 0 for $W_{a,j}^{\varepsilon}$, we obtain

$$\frac{m}{n} \sum_{t=1}^{[n/2]} E(\mu_{t,n} \hat{y}_{ab,t}^{\varepsilon} \hat{y}_{cd,t}^{\varepsilon})$$

$$\begin{aligned}
&= \frac{m}{n(w^*)^2}\left(\frac{1}{2\pi}\right)^2 \sum_{t=1}^{[n/2]} \sum_{j=-m/2}^{m/2} w_j^2 \mu_{t,n} \Phi_a(\exp(\lambda_{t+j})) \Phi_d'(\exp(-\lambda_{t+j})) \\
&\qquad\qquad\qquad\qquad \times \Phi_c(\exp(\lambda_{t+j})) \Phi_b'(\exp(-\lambda_{t+j})) + O(m/n) \\
&= \frac{m}{n(w^*)^2} \sum_{t=1}^{[n/2]} \sum_{j=-m/2}^{m/2} w_j^2 \mu_{t,n} f_{ad,t+j} f_{cb,t+j} + O(m/n) \\
&= \frac{m}{n(w^*)^2} \sum_{j=-m/2}^{m/2} w_j^2 \sum_{t=1}^{[n/2]} \mu_{t,n} f_{ad,t} f_{cb,t} + O(m/n).
\end{aligned}$$

The last equality is given by the Taylor expansion. Then, the assertion follows from Exercise 1.7.14 of Brillinger [2]. □

LEMMA 10. *Let* $\mu_i(\lambda), i = 1, 2,$ *be differentiable functions in* $[0, \pi]$ *and set* $\mu_{i,t,n} = \mu_i(\lambda_t).$ *Then, under* (A1)–(A4),

$$\begin{aligned}
\lim_{n\to\infty} \frac{m}{n} &\mathrm{Cov}\left(\sum_{t=1}^{[n/2]} \mu_{1,t,n} \hat{y}_{a_1b_1,t}^\varepsilon \hat{y}_{c_1d_1,t}^\varepsilon, \sum_{s=1}^{[n/2]} \mu_{2,s,n} \hat{y}_{a_2b_2,s}^\varepsilon \hat{y}_{c_2d_2,s}^\varepsilon \right) \\
&= \frac{D_u}{\pi} \int_0^\pi \mu_1(\lambda) \overline{\mu_2(\lambda)} \\
&\qquad \times (f_{a_1a_2}(\lambda) f_{b_2b_1}(\lambda) f_{c_1c_2}(\lambda) f_{d_2d_1}(\lambda) \\
&\qquad\quad + f_{a_1c_2}(\lambda) f_{d_2b_1}(\lambda) f_{c_1a_2}(\lambda) f_{b_2d_1}(\lambda))\, d\lambda
\end{aligned}$$

*for* $a_i, b_i, c_i, d_i = 1, \ldots, r, i = 1, 2,$ *where*

$$D_u = \frac{1}{2} \int_{-1/2}^{1/2} dx \int_{-1/2}^{1/2} dy \int_{-1}^1 u(x) u(y) u(x+z) u(y+z)\, dz \Big/ \left(\int_{-1/2}^{1/2} u(x)\, dx\right)^4.$$

PROOF. First, note that $\hat{y}_{ab,t}^\varepsilon$ and $\hat{y}_{cd,s}^\varepsilon$ are mutually independent for $|t-s| > m$, because $I_j^\varepsilon, j = 0, 1, \ldots, [n/2],$ are mutually independent variables. Then, noting that $E(\hat{y}_{ab,t}^\varepsilon) = 0$ and applying Lemmas 1 and 5 and Theorems 2.3.1 and 2.3.2 of Brillinger [2], we obtain

$$\begin{aligned}
\frac{m}{n} &\mathrm{Cov}\left(\sum_t \mu_{1,t,n} \hat{y}_{a_1b_1,t}^\varepsilon \hat{y}_{c_1d_1,t}^\varepsilon, \sum_s \mu_{2,s,n} \hat{y}_{a_2b_2,s}^\varepsilon \hat{y}_{c_2d_2,s}^\varepsilon \right) \\
&= \frac{m}{n} \sum_t \sum_s \mathrm{cum}(\mu_{1,t,n} \hat{y}_{a_1b_1,t}^\varepsilon \hat{y}_{c_1d_1,t}^\varepsilon, \overline{\mu_{2,s,n}} \hat{y}_{b_2a_2,s}^\varepsilon \hat{y}_{d_2c_2,s}^\varepsilon) \\
&= \frac{m}{n} \sum_t \sum_{s\,:\,|t-s|\le m} \mu_{1,t,n} \overline{\mu_{2,s,n}} [\mathrm{cum}(\hat{y}_{a_1b_1,t}^\varepsilon, \hat{y}_{c_1d_1,t}^\varepsilon, \hat{y}_{b_2a_2,s}^\varepsilon, \hat{y}_{d_2c_2,s}^\varepsilon)
\end{aligned}$$

$$+ \operatorname{cum}(\hat{y}^{\varepsilon}_{a_1b_1,t}, \hat{y}^{\varepsilon}_{b_2a_2,s}) \operatorname{cum}(\hat{y}^{\varepsilon}_{c_1d_1,t}, \hat{y}^{\varepsilon}_{d_2c_2,s})$$

$$+ \operatorname{cum}(\hat{y}^{\varepsilon}_{a_1b_1,t}, \hat{y}^{\varepsilon}_{d_2c_2,s}) \operatorname{cum}(\hat{y}^{\varepsilon}_{c_1d_1,t}\hat{y}^{\varepsilon}_{b_2a_2,s})]$$

$$= \frac{m}{n(w^*)^4} \sum_t \sum_{s,|t-s| \le m} \sum_{j_1,j_2,j_3,j_4} w_{j_1} w_{j_2} w_{j_3} w_{j_4} \mu_{1,t,n} \overline{\mu_{2,s,n}}$$

$$\times (f_{a_1a_2,t} f_{b_2b_1,t} f_{c_1c_2,t} f_{d_2d_1,t}$$

$$\times \delta_{t+j_1,s+j_3} \delta_{t+j_2,s+j_4}$$

$$+ f_{a_1c_2,t} f_{d_2b_1,t} f_{c_1a_2,t} f_{b_2d_1,t}$$

$$\times \delta_{t+j_1,s+j_4} \delta_{t+j_2,s+j_3}) + o(1)$$

$$= \frac{m}{n(w^*)^4} \sum_t \sum_{|t-s| \le m} \sum_{j_1,j_2} w_{j_1} w_{j_2} w_{j_1+t-s} w_{j_2+t-s} \mu_{1,t,n} \overline{\mu_{2,t,n}}$$

$$\times (f_{a_1a_2,t} f_{b_2b_1,t} f_{c_1c_2,t} f_{d_2d_1,t}$$

$$+ f_{a_1c_2,t} f_{d_2b_1,t} f_{c_1a_2,t} f_{b_2d_1,t}) + o(1).$$

The last two equalities follow from $2\pi j_i/n, i = 1, \dots, 4$ and $2\pi(t-s)/n$ being of order $m/n$, and, in the last summation, $j_1, j_2, s$ and $t$ must satisfy $|j_1 + t - s| \le m/2$ and $|j_2 + t - s| \le m/2$. Then, the assertion is obtained by Exercise 1.7.14 of Brillinger [2]. □

LEMMA 11. *Set*

$$\tau_{abcd,n} = \left(\frac{m}{n}\right)^{1/2} \sum_{t=1}^{[n/2]} c_{abcd,t,n} \hat{y}^{\varepsilon}_{ab,t} \hat{y}^{\varepsilon}_{cd,t}$$

*for* $a, b, c, d = 1, \dots, r$, *where* $c_{abcd,t,n}$ *are the constants bounded in* $t$ *and* $n$. *Then, under* (A1)–(A4),

$$\operatorname{cum}(\tau_{a_{11}b_{11}a_{12}b_{12},n}, \dots, \tau_{a_{k1}b_{k1}a_{k2}b_{k2},n}) = o(1)$$

*for* $k \ge 3$ *and* $1 \le a_{ij}, b_{ij} \le r, i = 1, \dots, k, j = 1, 2$.

PROOF. Applying Theorems 2.3.1 and 2.3.2 of Brillinger [2], we have

$$\operatorname{cum}(\tau_{a_{11}b_{11}a_{12}b_{12},n}, \dots, \tau_{a_{k1}b_{k1}a_{k2}b_{k2},n})$$

$$= \left(\frac{m}{n}\right)^{k/2} \sum_{t_1,t_2,\dots,t_k=1}^{[n/2]} \prod_{i=1}^{k} c_{a_{i1}b_{i1}a_{i2}b_{i2},t_i,n}$$

$$(21) \qquad \times \operatorname{cum}(\hat{y}^{\varepsilon}_{a_{11}b_{11},t_1} \hat{y}^{\varepsilon}_{a_{12}b_{12},t_1}, \dots, \hat{y}^{\varepsilon}_{a_{k1}b_{k1},t_k} \hat{y}^{\varepsilon}_{a_{k2}b_{k2},t_k})$$

$$= \left(\frac{m}{n}\right)^{k/2} \sum_{t_1,t_2,\dots,t_k=1}^{[n/2]} \prod_{i=1}^{k} c_{a_{i1}b_{i1}a_{i2}b_{i2},t_i,n}$$

$$\times \sum_{\nu^*} \prod_{l=1}^{p} \mathrm{cum}(\hat{y}_{a_{ij}b_{ij}, t_i}^{\varepsilon}; ij \in \nu_l),$$

where $\sum_{\nu^*}$ is over all the indecomposable partitions $\nu^* = \nu_1 \cup \cdots \cup \nu_p$ of the table

$$
\begin{array}{cc}
11 & 12 \\
21 & 22 \\
\vdots & \vdots \\
k1 & k2.
\end{array}
$$

Let $l^*$ be the number of elements of $\nu_l, l = 1, \ldots, p$. Since $E(\hat{y}_{ab,t}^{\varepsilon}) = 0$, any $l^*$ is greater than 1, and, hence, $p$ satisfies $p \le k$.

Now, applying Lemmas 4 and 5 to each cumulant on the right-hand side of (21) and noting that $\sum_{l=1}^{p} l^* = 2k$, we have

$$\prod_{l=1}^{p} \mathrm{cum}(\hat{y}_{a_{ij}b_{ij}, t_i}^{\varepsilon}, ij \in \nu_l) = O\left(\prod_{l=1}^{p} m^{1-l^*}\right) = O(m^{p-2k}).$$

By Theorem 2.3.1(iii) of Brillinger [2], the number of the nonzero terms of the summation $\sum_{t_1, t_2, \ldots, t_k}$ is $O(nm^{k-1})$, since $\hat{y}_{a_{ij}b_{ij}, t_i}^{\varepsilon}$ and $\hat{y}_{a_{i'j'}b_{i'j'}, t_{i'}}^{\varepsilon}$ are mutually independent for $|t_i - t_{i'}| > m$. Thus, the term on the right-hand side of (21) is $O(m^{p-k/2-1}/n^{k/2-1})$, which is $O((m/n)^{k/2-1}) = o(1)$ for $k \ge 3$. Hence, the proof is complete. $\square$

**Acknowledgments.** The authors are grateful to an Associate Editor and the referees for the valuable comments that helped us improve the manuscript significantly.

GRADUATE SCHOOL OF ECONOMICS
UNIVERSITY OF TOKYO
HONGO 7-3-1
BUNKYO-KU, TOKYO 113-0033
JAPAN
E-MAIL: yajma@e.u-tokyo.ac.jp

GRADUATE SCHOOL OF ECONOMICS
AND MANAGEMENT
TOHOKU UNIVERSITY
KAWAUCHI 27-1
AOBA-KU, SENDAI 980-8576
JAPAN
E-MAIL: matsuda@econ.tohoku.ac.jp